\documentclass{amsart}
\usepackage{a4wide,amssymb,color,hyperref}
\newcommand{\IZ}{\mathbb Z}

\theoremstyle{plain}

\theoremstyle{definition}

\theoremstyle{remark}

\newtheorem{example}{Example}[section]

\title{Steiner systems S(2,6,121/126), S(2,7,169) based on difference families}
\author{Ivan Hetman}
\date{January 2024}

\begin{document}

\begin{abstract} In this paper new Steiner systems $S(2,6,121)$, $S(2,6,126)$, $S(2,7,169)$ are introduced. Also some non-existence results for line lengths $7..11$ are presented. There is no solid proof that presented algorithm is exhaustive or correct, but it produces same results on already known difference families for line lengths $3..6$. Due to calculation-based approach this paper probably won't be published, but will be submitted to arxiv as it contains some new results.
\end{abstract}
\maketitle

\section{Introduction}
There exist enumerations for non-isomorphic difference families in \cite{HoCD} and \cite{BaiTop}. According to this results, number of non-isomorphic $(v,k,1)$ difference families are known for $k=3$ and $v\leq 61$, $k=4$ and $v\leq 76$, $k=5$ and $v\leq 85$, $k=6$ and $v\leq 96$ and $k=7$ and $v\leq 91$. Every difference family produces balanced incomplete block design $BIBD(v,k,\lambda)$. Number of BIBD is significantly larger than the family of difference families with same parameters. Still, in \cite{HoCD} only 1 non-isomorphic $BIBD(121,6,1)$ is known and only 2 non-isomorphic $BIBD(126,6,1)$.
\section{Algorithm description}
We will construct base blocks over $\IZ_{v}$ one by one. Blocks are lexicographically ordered. If $k|v$ then last block is always $(0,k,2k,...,(v - k))$. First item of block is always $0$, second item of block should produce maximal difference between two consecutive members of the block (including last and first = 0). Then on from second to last block we iterate limited possible numbers. Also we have filter which contains "forbidden" numbers which will make construction of difference family impossible. If specific block produces same difference $\Delta B = \{b_i - b_j | i, j = 1, ..., k; i \not = j\}$ as already computed, we filter it out because it will just repeat same calculation. In the end we add filtered difference family to the final list by "mirroring" procedure.

Difference family search algorithm is easily parallelizable by second number in first block. This fact helps to utilize modern multi-core hardware. Search was run on 24-core PC. Simple parallelization of the algorithm led to 10x performance boost. Code is present \cite{code} here.

Let's consider $v = 15$ and $k = 3$. As $3|15$, then we automatically have the last block $(0,5,10)$. Let's calculate first two blocks. After the search using base algorithm, we obtain first block $(0,7,9)$. Second block is then $(0,11,12)$. We see that this two blocks form difference family with $(0,5,10)$. Therefore it's desired difference family. As we limited first step by selecting only unique differences, we need to conduct "mirroring" procedure. To mirror any block, we fix first two numbers $b_1 = 0$ and $b_2$, and apply transformation $b_i \rightarrow b_1 + (v - b_i)$ to the rest. For instance, "mirror" of $(0,7,9)$ is $(0,7,13)$ and "mirror" of $(0,11,12)$ is $(0,11,14)$. If we mirror all blocks simultaneously it can be easily seen that obtained difference families are isomorphic. So, we need to conduct $2^{\lfloor\frac{v}{k(k-1)}\rfloor - 1}$ different mirrorings to obtain all possibly non-isomorphic difference families. In the end we filter out isomorphic difference families. For $v=15$ and $k=3$ nothing is filtered, but for $v=13$ and $k=3$, two difference families obtained by mirroring are isomorphic.

\section{Results}
Applying described algorithm gives us same numbers that are in \cite{HoCD} for known difference families. Still, number of non-isomorphic difference families for $v=73$ and $k=4$ in \cite{BaiTop} is $1426986$. But same calculation using algorithm described above gives us $1428546 > 1426986$ number.

Also using same algorithm it was possible to get some other low-hanging fruits. In addition to \cite{BaiTop1} result of non-existence of $(127,7,1)$ difference family, it was established non-existence of $(133, 7, 1)$, $(113, 8, 1)$, $(120, 8, 1)$, $(145, 9, 1)$, $(153, 9, 1)$, (Update December 2024) $(169, 8, 1)$, (Update December 2024) $(176, 8, 1)$, (Update December 2024) $(181, 10, 1)$, (Update December 2024) $(190, 10, 1)$ and (Update December 2024) $(231, 11, 1)$, (Update March 2025) $(175,7,1)$ difference families.

\begin{example} Exhaustive list of difference families $(121,6,1)$. "Mirrorings" of any of the sets produce 8 non-isomorphic difference family yielding 48 multiplier-nonisomorphic difference families.
\begin{itemize}
    \item \{\{0, 25, 37, 55, 76, 99\}, \{0, 52, 57, 72, 110, 113\}, \{0, 54, 71, 81, 90, 97\}, \{0, 73, 75, 79, 107, 108\}\}
    \item \{\{0, 25, 45, 66, 79, 97\}, \{0, 39, 46, 51, 83, 99\}, \{0, 62, 63, 65, 71, 98\}, \{0, 64, 74, 78, 93, 104\}\}
    \item \{\{0, 26, 35, 54, 72, 97\}, \{0, 40, 47, 52, 91, 113\}, \{0, 42, 53, 57, 63, 80\}, \{0, 56, 76, 89, 90, 92\}\}
    \item \{\{0, 26, 41, 53, 73, 96\}, \{0, 50, 52, 87, 111, 117\}, \{0, 63, 72, 77, 105, 108\}, \{0, 64, 75, 82, 83, 104\}\}
    \item \{\{0, 29, 31, 56, 80, 93\}, \{0, 40, 46, 50, 66, 98\}, \{0, 45, 60, 67, 78, 79\}, \{0, 47, 68, 77, 82, 85\}\}
    \item \{\{0, 38, 39, 56, 61, 85\}, \{0, 45, 51, 58, 66, 78\}, \{0, 49, 53, 79, 90, 93\}, \{0, 52, 54, 86, 102, 111\}\}
\end{itemize}
\end{example}
\begin{example} Exhaustive list of difference families $(126,6,1)$. "Mirrorings" of any of the sets produce 8 non-isomorphic difference family yielding 64 multiplier-nonisomorphic difference families. Last block $\{0,21,42,63,84,105\}$ is omitted.
\begin{itemize}
    \item \{\{0, 26, 38, 56, 81, 103\}, \{0, 40, 48, 68, 99, 102\}, \{0, 41, 57, 93, 94, 107\}, \{0, 80, 82, 87, 91, 97\}\}
    \item \{\{0, 28, 46, 68, 85, 101\}, \{0, 45, 52, 65, 95, 99\}, \{0, 49, 59, 64, 78, 115\}, \{0, 82, 88, 90, 91, 114\}\}
    \item \{\{0, 30, 50, 69, 95, 123\}, \{0, 34, 59, 77, 82, 94\}, \{0, 64, 68, 79, 119, 120\}, \{0, 80, 88, 90, 104, 117\}\}
    \item \{\{0, 35, 36, 64, 89, 103\}, \{0, 44, 60, 70, 77, 92\}, \{0, 55, 57, 61, 102, 107\}, \{0, 75, 83, 86, 95, 113\}\}
    \item \{\{0, 36, 37, 47, 77, 108\}, \{0, 38, 52, 64, 91, 97\}, \{0, 46, 50, 66, 69, 94\}, \{0, 68, 70, 75, 83, 92\}\}
    \item \{\{0, 36, 37, 64, 93, 115\}, \{0, 45, 52, 83, 106, 118\}, \{0, 67, 71, 77, 101, 117\}, \{0, 68, 82, 85, 87, 100\}\}
    \item \{\{0, 36, 38, 48, 75, 95\}, \{0, 52, 68, 81, 82, 86\}, \{0, 53, 62, 70, 77, 103\}, \{0, 55, 61, 80, 83, 115\}\}
    \item \{\{0, 47, 48, 67, 91, 123\}, \{0, 49, 66, 74, 111, 120\}, \{0, 57, 61, 68, 90, 95\}, \{0, 73, 85, 87, 103, 113\}\}
\end{itemize}
\end{example}

\begin{example} By exhaustive search finished in April 2025, I managed to find all difference families for $(169,7,1)$. All of them are obtained from one base family by mirrorings. Two of them have 3 multiplier automorphisms. Here is list of them together with multiplier automorphism count.
$$1: \{\{0, 1, 3, 11, 48, 65, 83\}, \{0, 4, 13, 29, 43, 81, 141\}, \{0, 5, 12, 36, 56, 78, 111\}, \{0, 6, 21, 40, 67, 90, 116\}\}$$
$$1: \{\{0, 1, 3, 11, 48, 65, 83\}, \{0, 4, 13, 29, 43, 81, 141\}, \{0, 5, 12, 36, 56, 78, 111\}, \{0, 6, 59, 85, 108, 135, 154\}\}$$
$$3: \{\{0, 1, 3, 11, 48, 65, 83\}, \{0, 4, 13, 29, 43, 81, 141\}, \{0, 5, 63, 96, 118, 138, 162\}, \{0, 6, 21, 40, 67, 90, 116\}\}$$
$$3: \{\{0, 1, 3, 11, 48, 65, 83\}, \{0, 4, 32, 92, 130, 144, 160\}, \{0, 5, 63, 96, 118, 138, 162\}, \{0, 6, 21, 40, 67, 90, 116\}\}$$
\end{example}

\section{Acknowledgement}

I would like to express my thanks to Roman Obukhivskyi who suggested better implementation for one of time-consuming algorithm parts, so it increased my algorithm speed 3-5 times and allowed to finish $(175,7,1)$ and $(169,7,1)$ on my home PC in acceptable amount of time (approximately a month).

\end{document}